\documentclass[11pt]{amsart}
\usepackage{mathrsfs}
\usepackage{amssymb,latexsym}

\usepackage{amssymb,amsmath,ams}
\newdimen\AAdi%
\newbox\AAbo%
\def\AAk#1#2{\s_etbox\AAbo=\hbox{#2}\AAdi=\wd\AAbo\kern#1\AAdi{}}%
\def\AAr#1#2#3{\s_etbox\AAbo=\hbox{#2}\AAdi=\ht\AAbo\raise#1\AAdi\hbox{#3}}%
\font\tenmsb=msbm10 at 11pt \font\sevenmsb=msbm7 at 8pt
\font\fivemsb=msbm5 at 6pt
\newfam\msbfam
\textfont\msbfam=\tenmsb \scriptfont\msbfam=\sevenmsb
\scriptscriptfont\msbfam=\fivemsb

\newlength{\abstractwidth}
\newcommand{\be}{\begin{equation}}
\newcommand{\bea}{\begin{eqnarray}}
\newcommand{\eea}{\end{eqnarray}}
\newcommand{\ee}{\end{equation}}

\def\ba{\begin{eqnarray}}
\def\ea{\end{eqnarray}}

\def\det{{\rm det}}

\def\14{{1 \over 4}}

\def\p{\partial}
\def\n{\nabla}

\def\e{\varepsilon}

\def\o{\omega}

\def\vp{\varphi}

\def\ddb{\partial\bar\partial}
\newtheorem{lemma}{Lemma}

\newtheorem{theorem}{Theorem}

\begin{document}

\title{PARTIAL LEGENDRE TRANSFORMS OF
NON-LINEAR EQUATIONS}
\author{Pengfei Guan}
\address{Department of Mathematics\\
         McGill University\\
         Montreal, Quebec. H3A 2K6, Canada.}
\email{guan@math.mcgill.ca}
\author{D. H. Phong}
\address{Department of Mathematics\\
Columbia University\\
New York, NY 10027, USA.}
\email{phong@math.columbia.edu }

\thanks{Research of the first author was supported in part
by NSERC Discovery Grant. Research of the second author was
supported in part by the National Science Foundation
grant DMS-07-57372.}

\begin{abstract}
The partial Legendre transform of a non-linear elliptic differential equation is shown
to be another non-linear elliptic differential equation. In particular,
the partial Legendre transform of the Monge-Amp\`ere equation
is another equation of Monge-Amp\`ere type. In
1+1 dimensions, this can be applied to obtain
uniform estimates to all orders for the degenerate
Monge-Amp\`ere equation with boundary data
satisfying a strict convexity condition.
\end{abstract}

\maketitle

\section{Introduction}
\setcounter{equation}{0}

The maximum rank and regularity properties of degenerate fully non-linear
equations are still largely unexplored, despite their considerable interest
for many geometric problems. For example, it is still an unresolved problem
raised by Donaldson \cite{D1, D2} to determine the precise regularity
of geodesics in the spaces of K\"ahler potentials and of volume forms
on a Riemannian manifold. These are given respectively by solutions of a degenerate
complex Monge-Amp\`ere equation and an equation introduced by Donaldson
\cite{D2}.
Many existence and regularity properties
have now been established for these
equations (see e.g. \cite{C00, PS1, CT, CH, H,
PS07, B, PS09, PS2} and references therein), but it is not known
how close they are to optimal.

\smallskip
In \cite{GP}, the maximum rank property has been established for several
special cases of the degenerate real Monge-Amp\`ere and Donaldson equations, for Dirichlet
data satisfying strict convexity conditions. Thus the natural question arises
of whether the regularity of solutions of fully non-linear elliptic
equations can be established, assuming that it is already known that they
have maximum rank.

\smallskip
A potentially useful feature of the maximum rank property is that it allows the use of a partial
Legendre transform. In fact, the partial Legendre transform was already exploited by
D. Guan \cite{Gd}
in determining geodesics for the space of K\"ahler potentials on toric varieties,
and by P. Guan \cite{Gp} and Rios, Sawyer, and Wheeden \cite{RSW1, RSW2}
in their study of the local regularity of certain degenerate Monge-Amp\`ere equations.
A first goal of this paper is to refine their analysis and show that,
even though the partial Legendre transform $f$ of a function $u$ is not a local
expression of $u$, the partial Legendre transform of an elliptic PDE in $u$
is another PDE in $f$ which is again elliptic. In particular, the
original Monge-Amp\`ere equation can be transformed globally into another {\it dual}
equation, again of Monge-Amp\`ere type, but which does not seem to have been
encountered before in the literature and may be of independent interest
(see Theorem \ref{1}, (b)). The partial strict convexity properties of one equation
are then equivalent to $C^2$ estimates for its dual, and one can expect a more general
correspondence between bounds for their derivatives. The second goal of this
paper is to apply this principle in the simplest case of the $1+1$ real Monge-Amp\`ere
equation (which coincides with the $1+1$ Donaldson equation). As a consequence,
we obtain $C^\infty$ bounds for this equation which depend only on the Dirichlet data, and in particular which remain uniform as the equation degenerates
(c.f. Theorem \ref{2}). We note that in this case, by \cite{Gd},
the solution of the limiting equation is already known to be smooth,
so the real interest of the result lies in the uniform validity of the approximation.

\section{Legendre transforms}
\setcounter{equation}{0}

The main goal of this section
is to work out the partial Legendre transforms of fully non-linear elliptic PDE's
in some generality. We shall find that, just as in the case of the full
Legendre transform, they are given by elliptic PDE's.

\subsection{The full Legendre transform}

We begin by re-visiting the standard Legendre
transform.
Let $u(x)$ be a strictly convex function on ${\bf R}^n$. Then $u(x)$ defines
a Legendre change of variables
\bea
\label{transform1}
x\to y=\frac{\p u}{\p x}(x).
\eea
Clearly, the Jacobian of this change of variables is
${\p y_j\over \p x_k}={\p^2 u\over \p x_j\p x_k}\equiv u_{jk}$.
The strict convexity implies that the map $x\to y$ from ${\bf R}^n$
to its image is invertible.
The Jacobian of the inverse $y\to x$ is given by the inverse $u^{jk}$ of
the Jacobian $u_{jk}$ of $x\to y$. Associated to the function $u(x)$
is also its Legendre transform $f(y)$, defined by
\bea
\label{transform0}
f(y)=xy-u(x),\qquad y={\p u\over\p x}.
\eea
Differentiating this relation with respect to $y$ shows that
the Legendre change of variables defined by $f(y)$ is the inverse map $y\to x$
\bea
\label{transform2}
y\to x={\p f\over \p y},
\eea
and we have the following exact analogues of the earlier formulas for $u$,
${\p y_j\over \p x_k}={\p^2 f\over \p x_j\p x_k}=u^{jk}$.
A partial differential equation of the form $F(u_{jk})=0$ can be viewed
as a partial differential equation in $f$. Its linearization has principal
symbol ${\p F\over \p u_{jk}}=f^{jp}f^{kq}\xi_p\xi_q$. Thus the ellipticity
of the equation in $u$ implies the ellipticity of the equation in $f$.
In particular, a Monge-Amp\`ere equation for $u$,
\bea
{\rm \det}\,({\p^2 u\over \p x_j\p x_k})=K
\eea
is equivalent to a Monge-Amp\`ere equation for $f$
\bea
{\rm \det}\,({\p^2 f\over \p y_j\p y_k})=K^{-1}.
\eea
We note that the changes of variables $x\to y$ and $y\to x$
in (\ref{transform1}, \ref{transform2}) are unaffected
if $u$ and/or $f$ are shifted by independent constants.
The relation (\ref{transform0}) can be viewed as a canonical way
of fixing the relative normalization of $f$ and $u$.

\subsection{The partial Legendre transform}

We come now to the situation of main interest to us,
namely partial Legendre transforms
of functions $u(x,t)$ which are periodic and satisfy a
a strict convexity condition in $x$.
More specifically,
let $e_i$ be the basis vectors for ${\bf R}^n$, that is, $e_i$ has component $1$ in the $i$-th position,
and all its other components are $0$. We consider functions $u(x,t)$ on
${\bf R}^n\times I$, $I=(0,1)$,
satisfying the periodicity condition
\bea
\label{shifts1}
u(x+ e_i)=u(x),
\qquad x\in {\bf R}^n,\quad 1\leq i\leq n,
\eea
and the strict convexity condition
\bea
\label{convexity}
{\p^2 u\over\p x_j\p x_k}+\delta_{jk}>0.
\eea
Thus $u$ can also be viewed as a function on $X\times I$, where $X$ is the torus
$X=({\bf R}/{\bf Z})^n$.

\medskip
We define the following Legendre change of variables
\bea
x\to y=(y_k),
\qquad y_k={\p u\over\p x_k}+x_k.
\eea
The inverse map $y\to x$ is well-defined and unique as a map from
${\bf R}^n$ to ${\bf R}^n$.
To see this, set $v=u+{1\over 2}|x|^2$,
and note that $v\to\infty$ as $x\to\infty$. Thus $v$ admits a minimum.
After a translation if necessary, we may assume that this minimum is at $0$,
and in particular, that $\n v$ vanishes there. It suffices then to see
that, for any $y$, the equation
\bea
{\p v\over \p x_k}(x_t)=ty_k,
\qquad
0\leq t\leq 1
\eea
can be solved by the method of continuity: the Jacobian $(v_{jk})$
is invertible implying openness of the set of $t$'s for which the equation
is solvable. And since $|\nabla v|\to \infty$ as $|x|\to\infty$,
we also find that $|x_t|$ is bounded, and the set of such $t$'s is also closed.

The maps $x\to y$ and $y\to x$ satisfy the following transformation
laws,
\bea
\label{shifts2}
y_k(x+e_i)=y_k(x)+\delta_{ik},
\qquad
x_k(y+e_i)=x_k(y)+\delta_{ik}.
\eea
The condition for $y_k$ follows immediately from its definition and the fact that $u(x)$
is periodic. To establish the condition for $y$, just observe that
\bea
y_k+\delta_{ik}
=
{\p u\over\p x_k}(x)+x_k+\delta_{ik}
=
{\p u\over\p x_k}(x+e_i)+(x+ e_i)_k
\eea
and the assertion follows by the uniqueness of the inverse map $y\to x$.

Clearly, the Jacobian of the map $x\to y$ is given by
\bea
{\p y_k\over\p x_j}={\p^2 u\over \p x_k \p x_j}+\delta_{jk}\equiv
g_{ij}.
\eea
Consequently, the Jacobian of the inverse map $y\to x$ is given by
\bea
{\p x_j\over\p y_k}=g^{jk}.
\eea
All these expressions are periodic, and descend to equations on the torus $X$.

\medskip
So far, we have discussed the Legendre maps defined by the function $u$.
We now define the Legendre transform $f$ of $u$ itself by the following formula,
\bea
\label{transform3}
f(y)=-{1\over 2}|x-y|^2 - u(x),
\qquad
y_j={\p u\over \p x_j}+x_j
\eea
for $y\in {\bf R}^n$. We note that, in view of the transformation laws
(\ref{shifts2}) for $x$ and $y$ under period shifts,
the function $f(y)$ is actually periodic,
and thus can be identified with a function on the torus $X$.
Again, the inverse map $y\to x$ can be viewed as the map associated with
the function $f$,
\bea
y\to x_j={\p f\over \p y_j}+y_j
\eea
and its Jacobian is given by
\bea
{\p x_k\over\p y_j}={\p^2 f\over\p y_j\p y_k}+\delta_{jk}=g^{jk}\equiv h_{jk}.
\eea
In particular, the Legendre transform $f$ satisfies the strict convexity
condition
\bea
D_y^2f+I>0
\eea
and we still have the relation
\bea
{\rm \det}\,({\p^2 u\over \p x_j\p x_k}+\delta_{jk})=
{\rm \det}^{-1}\,({\p^2 f\over \p y_j\p y_k}+\delta_{jk}).
\eea

Consider now the change in variables ${\bf R}^n\times I
\to {\bf R}^n\times I$ defined by
\bea
(x,t)\to (y,s),
\qquad y_j={\p u\over \p x_j}+x_j,
\ \ s=t.
\eea
The Jacobian of the inverse change of variables is given by
\bea
t_s=1,\qquad t_{y_p}=0,
\qquad
(x_k)_{y_p}=g^{kp},\qquad (x_k)_s=-g^{kp}u_{x_pt}.
\eea
It follows that the rule for differentiating a function $F(x,t)$ with respect to
the variables $(y,s)$ is given by
\bea
F_{y_p}=F_{x_j}g^{jp},
\qquad
F_s=-F_{x_j}g^{jp}u_{x_pt}+F_t.
\eea

The dependence of the partial Legendre transform on the additional
variables $s$ and $t$ is now conveniently described by the following three equations
\bea
\label{transform4}
\p_{y_j}u_t=-\p_s x_j,
\qquad
\p_s u_t=K\,({\rm \det}\,g)^{-1},
\qquad
f_s=-u_t.
\eea
Here $K$ is the $(n+1)\times (n+1)$ determinant ${\det}(D_{xt}^2 u+I_x)$,
\bea
K={\rm \det}\,(D_{xt}^2 u+I_x)=
u_{tt}({\rm \det} g)-G^{jp}u_{tx_j}u_{tx_p},
\eea
where $G^{jp}=({\rm det} g)g^{jp}$ is the matrix of co-factors of the metric
$g_{ij}$.
To see the first equation, we compute both sides. On the left, we
have $\p_{y_p}u_t=u_{tx_j}(x_j)_{y_p}=u_{tx_j}g^{jp}$. On the right,
we have $-\p_s x_p=-(-g^{pj}u_{tx_j})=g^{pj}u_{tx_j}$ also, as required.
Next, we apply the rule for differentiation of the previous paragraph and obtain
\bea
\p_s u_t= -u_{tx_j}g^{jp}u_{tx_p}+u_{tt}
=
({\rm \det} g)^{-1}(u_{tt}\,{\rm \det} g-G^{jp}u_{tx_j}u_{tx_p})
\eea
as claimed. Finally, differentiating the defining formula (\ref{transform3}) for $f$ gives
\begin{equation}
\p_sf=(y-x)\cdot x_s-(u_x \cdot x_s+u_t)
=
(y-x-u_x)\cdot x_s-u_t=-u_t.
\end{equation}
All three identities in (\ref{transform4}) have been proved.  They imply readily the following two
identities, which we also need later
\bea
\label{transform5}
u_{tx_j}=-f_{sy_k}h^{kj},
\qquad
u_{tt}=-f_{ss}+f_{sy_j}f_{sy_k}h^{jk}.
\eea

\subsection{The partial Legendre transform of non-linear
PDE's}

We consider now a fully non-linear equation of the form
\bea
F(D^2u)=0
\eea
on $X\times I$, where the unknown $u$ is required to satisfy the strict
convexity condition $D_x^2u+I_x>0$, and the equation is assumed to
be elliptic.
We would like to view this equation as an equation
for the Legendre transform $f$ of $u$. Note that $f$ is a non-local
quantity in $u$. Nevertheless, we have

\begin{theorem}
\label{1}
Let $X=({\bf R}/{\bf Z})^n$ be the $n$-dimensional torus,
and let $I=(0,1)$. Let $u(x,t)$ be a function
on $X\times I$ satisfying the strict convexity condition
(\ref{convexity}). Let
$(x,t)\to (s=t,y)$ be the partial Legendre transform as defined by
(\ref{transform2}), and let
$f$ be the partial Legendre transform of the function $u$
as defined by (\ref{transform3}).

\rm{(a)} If $F(D^2u)=0$ is a second-order elliptic PDE for $u$,
then it can also be viewed as a second-order elliptic PDE for
the partial Legendre transform $f$ of $u$.

\rm{(b)} If particular,
$u$ satisfies the Monge-Amp\`ere equation
\bea
\label{MAequation}
{\rm \det}\,(D_{x,t}^2 u+I_x)=K
\eea
if and only if its partial Legendre transform $f$
satisfies the following equation on $X\times I$
also of Monge-Amp\`ere type,
\bea
\label{main-eq}
{\p^2 f\over \p s^2}
+
K\, {\rm \det}\,(D^2_y f+I)=0.
\eea
\end{theorem}

\bigskip
{\noindent}
{\it Proof}:
From the discussion in the preceding section, the
partial Legendre transform $f$ of $u$ is a well-defined function on $X\times I$,
$\p_j\p_k u+\delta_{jk}$ are given by the inverse
of the matrix $\p_j\p_k f+\delta_{jk}$, and $u_{tx_j}$
and $u_{tt}$ are given by the expressions
(\ref{transform5}) in the second derivatives of $f$. Thus the equation $F(D^2u)$
is automatically a second order non-linear equation for $f$.
To verify the ellipticity of the equation viewed as an equation for $f$,
we work out first the linearized operator of $F(D^2u)$,
keeping variations in $\delta u$,
\bea
\label{linearizationF}
\delta F
={\p F\over\p u_{tt}}(\delta u)_{tt}
+
{\p F\over\p u_{tx_j}}(\delta u)_{tx_j}
+
{\p F\over\p u_{x_jx_k}}(\delta u)_{x_jx_k}.
\eea
We need to express this quantity in terms
of the derivatives of $\delta f$. In view
of the expression (\ref{transform5}) for $\delta u_{tt}$ and $\delta u_{tx_j}$,
we have
\bea
\delta u_{tt}&=&
-
\delta f_{ss}+2\delta f_{sy_j}f_{sy_k}h^{jk}-f_{sy_j}f_{sy_k}\delta f^{jk}
\nonumber\\
\delta u_{tx_j}
&=&
-\delta f_{sy_k} h^{kj}+f_{sy_j}\delta f^{jk},
\eea
where the indices are raised or lowered using the metric $h_{jk}$.
Thus we have
\bea
\delta F
&=&
-{\p F\over\p u_{tt}}\delta f_{ss}
+
(-{\p F\over \p u_{tx_j}}+2{\p F\over \p u_{tt}}f_{sy_j})h^{jk}\delta f_{sy_k}
\nonumber\\
&&
+
(-{\p F\over\p u_{tt}}f_{sy_j}f_{sy_k}
+
{\p F\over\p u_{tx_j}}f_{sy_k}
-
{\p F\over \p u_{x_jx_k}})\delta f^{jk}.
\eea
This means that,
if $\tau$ and $\xi_j$ are respectively
the variables dual to $t$ and $x_j$, the symbol $\sigma(\tau,\xi)$ of the linearized operator
is given by
\bea
\sigma(\tau,\xi)
&=&
{\p F\over\p u_{tt}}\tau^2
+
({\p F\over \p u_{tx_j}}-2{\p F\over \p u_{tt}}f_{sy_j})h^{jk}\tau\xi_k
\nonumber\\
&&
+
({\p F\over\p u_{tt}}f_{sy_j}f_{sy_k}
-
{\p F\over\p u_{tx_j}}f_{sy_k}
+
{\p F\over \p u_{x_jx_k}})\xi_j\xi_k.
\eea
We shall show that $\sigma(\tau,\xi)$ is positive definite
if the original equation $F(D^2u)=0$ is elliptic.
Introduce the variable $\eta_j\equiv h^{jk}\xi_k$ for convenience.
Then completing the square in $\tau$ gives
\bea
\label{symbol}
\sigma(\tau,\xi)
&=&
[\tau \sqrt{\p F\over\p u_{tt}}+
({1\over 2}{\p F\over\p u_{tx_j}}-{\p F\over\p u_{tt}}f_{sy_j}){\eta_j\over
\sqrt{\p F\over \p u_{tt}}}]^2
\nonumber\\
&&
+
{1\over 4{\p F\over\p u_{tt}}}\,
[4{\p F\over\p u_{tt}}\,{\p F\over\p u_{x_jx_k}}\eta_j\eta_k
-
({\p F\over\p u_{tx_j}}\eta_j)^2].
\eea
To prove part (a) of the Theorem, it suffices then to show that
that the second term on the right is strictly positive for $\eta\not=0$
when $F$ is elliptic. Since the symbol of the linearized operator
when the unknown is $u$ is given by
\bea
{\p F\over\p u_{tt}}\tau^2
+
{\p F\over \p u_{tx_j}}\tau\xi_j
+
{\p F\over\p u_{x_jx_k}}\xi_j\xi_k
\eea
its ellipticity does imply the positivity of the second term on the right
hand side of (\ref{symbol}). Part (a) is proved.

Part (b) follows immediately from the identities in (\ref{transform4}),
\bea
f_{ss}=-\p_su_t=-K({\rm \det} \,g)^{-1}
=-K\,{\rm \det}\,(D_y^2f+I).
\eea
The proof of Theorem \ref{1} is complete. Q.E.D.

\bigskip
We conclude this section with a few remarks.

\medskip

$\bullet$ Part (b) of Theorem \ref{1} can be viewed as a refinement of several earlier results in the literature using the partial Legendre transform: when $K=0$, it reproduces the result of D. Guan \cite{Gd}.
For general $K$, and when the considerations are local (instead of on a torus as
here), then P. Guan \cite{Gp} in dimension $n=1$
and Rios-Sawyer-Wheeden \cite{RSW1} for general dimension $n$
have shown that the coordinates $x_j$, viewed as functions
of $(y,s)$, satisfy the following elliptic,
non-linear system of equations
\bea
\p_s^2 x_j+{\p\over \p y_j}(K{\rm \det}({\p x_k\over\p y_m}))=0,
\qquad 1\leq j\leq n.
\eea
This system of equations follows immediately from differentiation of the equation (\ref{main-eq}) with respect to $y_j$.

$\bullet$ The presence of the background symmetric form $\delta_{jk}$
is very similar to the presence of the K\"ahler form $\o$ for the complex
Monge-Amp\`ere equation $(\o+{i\over 2}\ddb \phi)^n=F\o^n$.

$\bullet$ The correspondence between an equation in $u$ and
its ``dual" equation in $f$ can provide non-trivial information.
For example, it is not evident that the equation
(\ref{main-eq}) admits smooth solutions for given Dirichlet
data, even when $K$ is a strictly positive constant. On the
other hand, the existence of such solutions is an immediate consequence
of the existence of smooth solutions for the dual equation
(\ref{MAequation}), which can be established by the
theory of Caffarelli-Nirenberg-Spruck \cite{CNS}, with
the improved barrier arguments of B. Guan \cite{Gb}.

Similarly, lower bounds for $D_x^2u+I$ are equivalent to upper bounds
for $D_y^2f+I$ and vice versa, and the problems of partial $C^2$ estimates
and partial strict convexity are in this sense ``dual".
For example, the $C^2$ estimates for the original Monge-Amp\`ere equation
(\ref{MAequation})
can be established by traditional methods as in \cite{CNS},
or as a consequence of the convexity results for the dual equation
(\ref{main-eq}), using for example the recent results of Bian-Guan \cite{BG1, BG2}.

$\bullet$ Beyond the Monge-Amp\`ere equation,
the partial Legendre transforms of Hessian equations may be of interest. By (\ref{transform5}),
\bea
D_{x,t}^2u+I=\begin{pmatrix}-\tilde K({\rm \det}\,\tilde g)^{-1}
& -\tilde u_{sy_1}\lambda^{-1}_1&\cdots &-\tilde u_{sy_n}\lambda^{-1}_n\cr
-\tilde u_{sy_1}\lambda^{-1}_1&\lambda^{-1}_1& \cdots &0 \cr
\cdots \cr
-\tilde u_{sy_n}\lambda^{-1}_n&0 &\cdots &\lambda^{-1}_n\cr\end{pmatrix}
\eea
where $\lambda_i$ are the eigenvalues of $\tilde g=(\tilde u_{ij}+\delta_{ij})$,
and we have denoted for convenience all quantities associated
with the partial Legendre transform by a $\tilde{}$
(e.g. $f$ is now denoted $\tilde u$, and $\tilde K$
is the Monge-Amp\`ere determinant of $D_{y,s}^2\tilde u+I$).
The standard formula for the $k$-th symmetric function $\sigma_k$
of the eigenvalues of a matrix is
\begin{equation}
\sigma_k(V)={1\over k!}\sum \delta_{j_1\cdots j_k}^{i_1\cdots i_k}v_{i_1j_1}\cdots v_{i_kj_k}
\end{equation}
We apply this formula to the the symmetric function $\hat\sigma$ of the $n+1$-dimensional
matrix $D_{x,t}^2u+I$.
We find
\bea
\hat\sigma_k(D_{x,t}^2u+I)=\sigma_k(\tilde g^{-1})-
\sum_{i_1\cdots i_{k-1}}\lambda^{-1}_{i_1}\cdots\lambda^{-1}_{i_{k-1}}
(\sum_{i_j=1}^{k-1}\tilde u_{sy_{i_j}}^2\lambda^{-1}_{i_j}+\tilde K(\det \,\tilde g)^{-1}).
\eea
From here, it is easily seen that the Laplace equation
$u_{tt}+\Delta u=K_1$ gets transformed into
\begin{equation}
\tilde K +K_1 {\rm det}\,\tilde g =\sigma_{n-1}(\tilde g)
\end{equation}
For general $k$,
the equation $\hat\sigma_k(u)=K_k$ gets transformed into
\begin{equation}
K_k={\sigma_{n-k}(\tilde g)\over\sigma_n(\tilde g)}
-
\sum_{i_1\cdots i_{k-1}}\lambda^{-1}_{i_1}\cdots\lambda^{-1}_{i_{k-1}}
(\sum_{i_j=1}^{k-1}\tilde u_{sy_{i_j}}^2\lambda^{-1}_{i_j}+\tilde K(det \,\tilde g)^{-1})
\end{equation}
We note that when $k=(n+1)$, the above identity recovers (\ref{main-eq}).

\section{The $1+1$ Monge-Amp\`ere equation}
\setcounter{equation}{0}

In this section, we consider more specifically the case $n=1$ of the Monge-Amp\`ere
equation. In this case, the equation becomes
the following equation on $X\times I$,
\bea
\label{1+1-u}
u_{tt}(1+u_{xx})-u_{xt}^2=\e
\eea
where $X={\bf R}/{\bf Z}$ is a circle and $\e>0$ is a constant. The solution
$u$ is required to satisfy $D_{xt}^2u+I_x\geq 0$,
and we impose the Dirichlet condition $u(t,0)=u^0(x)$, $u(t,1)=u^1(x)$,
with $u^i\in C^\infty(X)$ satisfying the strict convexity condition,
\bea
\label{Dirichlet-u}
u^i_{xx}+1\geq \lambda
\eea
for $i=0,1$ and $\lambda>0$ a constant. We have

\begin{theorem}
\label{2}
Let $u$ be the solution of the equation (\ref{1+1-u}) on $X^1\times I$,
with $D_{xt}^2u+I_x\geq 0$ and smooth Dirichlet data satisfying
(\ref{Dirichlet-u}). Then for any non-negative integer $N$,
there exists $M(N)$ and a constant $C_N$ depending only
on  $\lambda>0$ and the $C^{M(N)}$ norms of the Dirichlet data $u^0$, $u^1$
so that
\bea
\sum_{a+b\leq N}\|\p_x^a\p_t^bu\|_{C^0(X^1\times I)}
\leq C_N.
\eea
In particular, the constants $C_N$ are independent of $\e$.
\end{theorem}

The rest of this section is devoted to the proof of Theorem \ref{2}.
To apply the partial Legendre transform as in \S 2, we need the strict
partial convexity of $u$. This follows from Theorem 1 of \cite{GP}.
But since the present case is particularly simple, we can supply the short
proof for the convenience of the reader:

\begin{lemma}
\label{maximumrank}
Let $u(x,t)$ be the solution of the equation (\ref{1+1-u}),
as specified in the statement of Theorem \ref{2}.
Then $u_{xx}(x,t)+1\geq\lambda$ for all $x\in X$.
\end{lemma}

\noindent
{\it Proof of Lemma \ref{maximumrank}}:  Let $\lambda_0
={\rm min}_{(x,t)\in X\times\bar I}(u_{xx}+1)$, and set
$\vp(x,t)=u_{xx}+1-\lambda_0$.
We establish a strong maximum principle for $\vp$.
If $\vp$ vanishes on the boundary, the lemma is proved.
We shall show that in a neighborhood of any interior zero of $\vp$,
$\vp$ satisfies an elliptic differential inequality equation of the form
\bea
\label{strongmaximum}
F^{ij}\vp_{ij}\leq C\,|\n\vp|.
\eea
Here we have denoted by $F$ the function
$F(D_{xt}^2u+I_x)$ with $F(M)={\rm det}(M_{ij})$,
for $M$ any symmetric $2\times 2$ matrix, and $F^{ij}={\p F\over \p M_{ij}}$.
The constant $C$ is required to be independent of the point $(x,t)$ but may depend on everything else.
By the strong maximum principle, this would imply that $\vp$ vanishes in a neighborhood of
any interior zero. If the set of such zeroes is not empty, then $\vp$ vanishes
identically. By continuity $\vp$ would again vanish on the boundary, and the lemma is proved
in all cases.

The equation (\ref{1+1-u}) can be written as
$F(D_{xt}^2u+I_x)=\e$. Differentiating the equation successively gives
\bea
F^{ij}u_{ijx}=0,
\qquad
F^{ij}u_{ijxx}+F^{ij,kl}u_{ijx}u_{klx}=0.
\eea
Thus $F^{ij}\vp_{ij}=-F^{ij,kl}u_{ijx}u_{klx}$, and more explicitly,
\bea
F^{ij}\vp_{ij}
=-2 u_{ttx}u_{xxx}-2 u_{xxt}^2=-2u_{ttx}\vp_x-2\vp_x^2.
\eea
The inequality (\ref{strongmaximum}) follows. Q.E.D.

\bigskip

Let $f(y,s)$ be now the partial Legendre transform
of $u$, as defined in section. When $n=1$, the equation for $f$ simplifies to
\bea\label{laplace-equation}
Lf\equiv f_{ss}+\e f_{yy}=0\quad{\rm on}\ \ X\times I,
\eea
and $f(y,0)$ and $f(y,1)$ are given by the Legendre tranforms $f^0(y)$ and $f^1(y)$
of the functions $u^0(x)$ and $u^1(x)$.
General linear elliptic theory says that any derivative
of $f$ can be bounded in terms of the boundary data and
the ellipticity constant $\e$. However, we require estimates which
are uniform as
$\e\to 0$, and such estimates do not
seem to have been written down in the literature.
We provide below a brief and explicit derivation of estimates
uniform in $\e$,
exploiting the simple form of the equation and of the boundary
in the present case.  More precisely, we shall establish the
following lemma:

\begin{lemma}
\label{laplacian}
Consider the Dirichlet problem for the Laplacian $L$ in
(\ref{laplace-equation}) with $\e$ a constant satisfying
$0<\e<1$. Then for any $m,k$ we have
\bea
\|\p_y^m\p_s^k f\|_{C^0(\tilde X\times I)}
\leq C_{m,k}
\eea
where $C_{m,k}$ are constants which depend only
on the Dirichlet data (and on $m$ and $k$). In particular, $C_{m,k}$ are
independent of $\e$.
\end{lemma}

\noindent
{\it Proof}:
Clearly $\p_y^m\p_s^k f$ satisfies the same Laplace equation. By the
maximum principle, we have then
\bea
\|\p_y^m\p_s^k f\|_{C^0(\tilde X\times I)}=
\|\p_y^m \p_s^kf\|_{C^0(\tilde X\times \p I)}.
\eea
We shall show the right hand sides can be estimated in terms of
the Dirichlet data alone, for arbitrary $m$ and $k=0$ or $k=1$.
Assuming this, it follows that the left hand side is also bounded in terms of
the Dirichlet data alone for these values of $m$ and $k$.
Since the equation implies that
$\p_y^m\p_s^k f=-\e\p_y^{m+2}\p_s^{k-2}f$
for $k\geq 2$, it follows that uniform bounds
for $\p_y^m\p_s^k f$ for arbitrary $k$ follow from
the special cases $k=0$ and $k=1$, and the lemma would
be proved.

\smallskip
We return to the proof of bounds for $\|\p_y^m\p_s^k f\|_{C^0(\tilde X\times\p I)}$
when $k=0$ or $k=1$. When $k=0$, they are obvious,
so we concentrate on the case $k=1$.
Let $\tilde f$ be a function with the same boundary values as $ f$
(e.g., $\tilde f= t f(y,1)+(1-t) f(y,0)$).
Set
\bea
w(y,s)=-As^2+Bs
\eea
for constants $A,B>0$ which we shall choose in a moment.
Since $L(\p_y^mf)=0$, we have
\bea
L(\p_y^m(f-\tilde f))
-L(\p_y^m\tilde f)
\geq c_1
\eea
where $c_1$ is a constant depending only on the Dirichlet data
and independent of $0<\e<1$. On the other hand
\bea
Lw=-2A
\eea
so we can choose $A$ large enough so that $Lw\leq L(\p_y^m(f-\tilde f))$
on $\tilde X\times I$. Now $\p_y^m(f-\tilde f)$
vanishes identically on the boundary, so if we choose
$B$ large enough so that $w\geq 0$ everywhere,
we shall also have $\p_y^m(f-\tilde f)\leq w$
on the boundary. Thus, by the comparison principle,
$\p_y^m(f-\tilde f)\leq w$ on $\tilde X\times I$.
Since both of these functions vanish at $s=0$, we obtain
\bea
\p_s\p_y^mf \leq -2As+B+\p_s\p_y^m\tilde f,
\eea
which is an upper bound for $\p_s\p_y^mf$. Applying the argument
to $-f$ instead of $f$, we obtain a lower bound for
$\p_s\p_y^mf$ at the boundary points $(y,0)$. Since the argument at the
boundary points $(y,1)$ is identical,
the proof of the lemma is complete.

\bigskip
Next, we show that $C^N$ bounds for $x$ (viewed as a function
of $(y,s)$, together with a strict partial
convexity bound, imply $C^M$ bounds for the original function $u$:

\begin{lemma}
\label{bounds f-u}
For any non-negative integer $M$, we have
\bea
\sum_{m+b\leq N}\|\p_x^m\p_t^bu\|_{C^0(X^1\times I)}
\leq C_M
\eea
where $C_M$ is a constant depending only on the Dirichlet data $u^0, u^1$,
the lower bound $\lambda>0$, and the $C^0$ norm of a finite number $N(M)$
of spatial derivatives $\p_y^m x$ of the function $x$.
\end{lemma}

{\it Proof of Lemma \ref{bounds f-u}}:
First, we show that bounds for $\p_y^m x$
imply bounds for $\p_x^a u$.
This is an easy consequence of the following formula,
which is itself a consequence of the chain rule
established in the previous section:
\bea
\label{higherderivatives}
\p_y^mx=
-{1\over (u_{xx}+1)^m}\p_x^{m+1}u
+
{P(u,\cdots,\p_x^m u)\over (u_{xx}+1)^{2m-1}}
\eea
where $P$ is a generic notation for
a polynomial in all its entries for all $m\geq 2$.

Now for $m=2$, the bounds of $\p_x^mu$
in terms of boundary data alone are a special case of
the $C^2$ estimates for the Dirichlet problem for the Monge-Amp\`ere
equation \cite{CNS}. Note that these bounds do not require a strictly
positive lower bound for the Monge-Amp\`ere determinant,
and thus give bounds which are uniform in $\e$ in our case.

Assume that bounds depending only on the Dirichlet data and
a strictly positive lower bound $\lambda$ for $u_{xx}+1$ have been established
for $m$. In view of the above formula for $\p_y^mx$,
it follows that such bounds for $\p_x^{m+1}u$
reduce to such bounds for $\p_y^m x$ on $\tilde X\times I$.
By the maximum principle for $\p_y^m x$, this reduces in turn
to bounds for $\p_y^m x$ only on the boundary.
But the same formula (\ref{higherderivatives}) for higher
derivatives above shows that
$\p_y^m x$ on the boundary $\tilde X\times\p I$ is
determined completely by the boundary data.
This establishes the desired bounds
for $\|\p_x^a u\|_{C^0(X\times I)}$.

\medskip
Next, we consider mixed derivatives
of the form $\p_x^a\p_tu$. It is again easy to establish
the following general formula linking $\p_y^m\p_s x$ and
$\p_x^a\p_t u$ for all $m\geq 0$,
\begin{equation}
\p_y^m\p_s x=
-{\p_x^{m+1}\p_t u\over (u_{xx}+1)^{m+1}}
+
{P(u,\cdots,\p_x^{m+2}u,\p_tu,\cdots,\p_x^m \p_tu)
\over (u_{xx}+1)^{2m+1}},
\end{equation}
where $P$ is again a polynomial in all its entries.
In view of Lemma \ref{laplacian}, the left hand side can be bounded
uniformly in terms of the Dirichlet data and the lower bound
$\lambda$ for $u_{xx}+1$.
Thus the formula implies that $\p_x^{m+1}\p_t u$ can be
similarly bounded if $\p_x^m\p_t u$ is. Since $u_{xt}$ is bounded by the Dirichlet
data in view of the $C^2$ estimates, we obtain by induction
the uniform boundedness of $\p_x^m\p_t u$ for all $m$.

\medskip
Finally, by differentiating the original Monge-Amp\`ere equation,
we can show inductively on the number $b$ of $t$ derivatives
in $\p_x^a\p_t^b u$ that they are in turn bounded.
First, we show this for $b=2$.
Differentiating the equation $m$ times with respect to
$x$ gives
\begin{equation}
\p_x^m u_{tt}=
{1\over u_{xx}+1}
P(\p_x^2u,\cdots,\p_x^{m+2}u,
\p_t\p_xu,\cdots,\p_t\p_x^{m+1}u,
\p_t^2u,\cdots,\p_t^2\p_x^{m-1}u)
\end{equation}
with $P$ a polynomial in all its entries.
Since $\p_x^m u_{tt}$ is bounded for $m=0$ by the $C^2$
estimates, and $\p_x^a u$, $\p_x^a\p_t u$ are now known to be bounded,
the formula allows us to show by induction on $m$
that $\p_x^m \p_t^2 u$ is bounded.
Next, assume that $\p_x^m\p_t^b u$ is bounded for an integer $b\geq 2$ and all
non-negative integers $m$.
We shall show that this remains true if $b$ is replaced by $b+1$. Differentiating the
equation $b-1$ times gives
\bea
\p_t^{b+1}u={1\over u_{xx}+1}
P(\p_x \p_t,\cdots,\p_x\p_t^b u,\p_x^2\p_t u,\cdots,\p_x^2\p_t^b u)
\eea
where $P$ is a polynomial in all its entries. This shows that $\p_x^m\p_t^{b+1}u$
is bounded for $m=0$. Thus it suffices to show that if $\p_x^\ell\p_t^{b+1}u$
is bounded for all non-negative
integers $\ell\leq m$, then the same is true
if $m$ is replaced by $m+1$. But this follows readily by differentiating
the previous formula,
\begin{equation}
\p_x^{m+1}\p_t^{b+1}u=\p_x^{m+1}
\big\{
{1\over u_{xx}+1}
P(\p_x \p_t,\cdots,\p_x\p_t^b u,\p_x^2\p_t u,\cdots,\p_x^2\p_t^b u)
\big\}.
\end{equation}
Since the right hand side involves terms with at most $b$ derivatives in $t$,
and all derivatives in $x$ of such expressions are bounded, the desired bound
for $\p_x^{m+1}\p_t^{b+1}u$ follows. Q.E.D.

\bigskip
Putting together Lemmas 2-4, we obtain
Theorem \ref{2}.

\bibliographystyle{amsplain}

\end{document}